\documentclass[12pt]{amsart}

\title[Random walks on the torus]{Random walks on the torus with several generators}
\author[Timothy Prescott]{Timothy Prescott$^*$}
\thanks{$^*$Department of Mathematics, University of 
California, Los Angeles, CA 90095, {\tt tmpresco@math.ucla.edu}.}
\author[Francis Edward Su]{Francis Edward Su$^{**}$}
\thanks{$^{**}$Research partially supported by NSF Grant
  DMS-0301129.  Department of Mathematics, Harvey Mudd College,
  Claremont, CA  91711, {\tt su@math.hmc.edu} (corresponding author)}

\newcommand{\R}{\mathbb{R}}
\newcommand{\T}{\mathbb{T}}
\newcommand{\Z}{\mathbb{Z}}

\newtheorem{theorem}{Theorem}
\newtheorem{lemma}[theorem]{Lemma}
\newtheorem{defn}[theorem]{Definition}

\begin{document}

\baselineskip=.27in

\begin{abstract}
Given $n$ vectors $\{\vec\alpha_i \}_{i=1}^n \in [0,1)^d$,
consider a random walk on the $d$-dimensional torus $\T^d = \R^d /
\Z^d$ generated by these vectors by successive addition and
subtraction.
For certain sets of vectors, this walk converges to Haar (uniform)
measure on the torus.
We show that the 
discrepancy distance $D(Q^{*k})$ between the $k$-th step distribution of the
walk and Haar measure is bounded below 
by $D(Q^{*k})\geq C_1 k^{-n/2}$, where $C_1=C(n,d)$ is a constant.
If the vectors are badly
approximated by rationals (in a sense we will define) then
$D(Q^{*k})\leq 
C_2 k^{-n/2d}$ for
$C_2=C(n,d,\vec\alpha_j)$ a constant.
\end{abstract}

\maketitle

Let $\T^d = \R^d / \Z^d$ denote the $d$-dimensional torus.  
As a quotient group of $\R^d$ it is an additive group, so the group
elements may be viewed as elements of $[0,1)^d$, with the group
operation defined as coordinate-wise addition mod 1.

Let $\vec\alpha_1, \vec\alpha_2, \ldots, \vec\alpha_n$ be vectors in
$\T^d = [0,1)^d$,
and consider the random walk on the $d$-dimensional 
torus $\T^d$ that proceeds as follows.  Start at $\vec 0$.
At each step, choose one the vectors $\vec\alpha_i$ with probability
$1/n$ and add or subtract that vector 
(with probability $1/2$)
to the current position to get to the next position in the walk.

As a random walk on a group, 
the $k$-th step distribution of the walk converges to a limiting
distribution~\cite{kloss}, and in many cases this will be Haar
measure, the unique translation-invariant measure on the group.  For
the torus $\T^d$, Haar measure may be thought of as the uniform distribution
on the ``flat'' cube $[0,1)^d$, since addition corresponds to translation on
$\R^d/\Z^d$.
We shall prove bounds for how quickly this random walk approaches Haar
measure on the torus.

We first note that for certain sets of vectors, this walk
may not converge to Haar measure.  For instance, if all the entries of
each $\vec\alpha_i$ are rational, then the random walk will not
converge to Haar measure, but will converge to a limiting distribution
supported on a discrete subgroup of $\T^d$.  As another example, if
there is only one generator $\vec\alpha_1=(x,x,...,x)$ for some irrational
$x$, then the walk will be supported on a circle along 
the ``diagonal'' of the torus.  (However, a single vector can
generate a walk that does converge to Haar measure, provided it is chosen
well.)

Let $Q$ denote the generating measure for this random walk, i.e., 
if $S =\cup_{i=1}^n \{+\vec\alpha_i,-\vec\alpha_i\}$ is the set of 
generators of the random walk, then 
for a set $B\subseteq\T^d$, let $Q(B)=| B\cap S|/|S|$ where $|\cdot|$
denotes the size of a finite set.
The $k$-th step probability distribution is then given by the $k$-th
convolution power of $Q$, which we denote by $Q^{*k}$.  Let $U$ denote
Haar measure.

As a measure of distance between the probability distributions
$Q^{*k}$ and $U$, we will use
the {\em discrepancy metric}, which is defined to be the supremum of
the difference of two probability measures over all ``boxes'' in
$\T^d=\R^d/\Z^d$ with sides parallel to the axes in $\R^d$,
i.e., of the form $[a_1,b_1)\times [a_2,b_2)\times ... \times [a_d,b_d)$.
Let $D(Q^{*k})$ denote the discrepancy of $Q^{*k}$ from Haar measure $U$:
$$
D(Q^{*k}) := \sup_{box B \subseteq \T^d} |Q^{*k}(B)-U(B)|.
$$
The discrepancy metric has been used by number theorists to study the 
uniform distribution of sequences mod 1, e.g.,
see~\cite{drmota-tichy,kuipers-nied}.  Diaconis~\cite{diaconis} suggested
its use for the study of rates of convergence for random walks on
groups.  
It admits Fourier bounds~\cite{hensley-su} and has many other
nice properties and connections with other probability
metrics~\cite{gibbs-su}.

Although the {\em total variation} metric is more commonly 
used to study the convergence of random walks, we do not
use it here because this random walk does not converge in total
variation (in fact, the total variation distance between $Q^{*k}$ and
$U$ is always 1, since at any step $Q^{*k}$ is supported on a finite
set).  The possibility of using Fourier analysis to bound the
discrepancy distance makes it a more desirable choice than 
other common metrics on probabilities, such as the Prohorov
metric, and has allowed many recent results for the study of discrete
random walks on continuous state spaces (e.g., \cite{su-tams,
  su-drunkard}).  Most of the literature for rates of convergence of
random walks have been limited to walks on finite groups or state spaces, and
those that have focused on infinite compact groups 
(e.g., \cite{porod}, \cite{rosenthal}, \cite{su-leveque}) 
have studied walks generated by
continuous measures.  By contrast, the walk we study is generated by a
discrete set of generators on an infinite group.  

We prove:

\begin{theorem}
\label{lowerbound}
Let $Q$ denote the generating measure of the 
the random walk on the $d$-torus generated by $n$ vectors 
$\vec\alpha_1,...,\vec\alpha_n$.  Then the $k$-th step
probability distribution $Q^{*k}$ satisfies:
$$
D(Q^{*k}) \ge 
\frac{1}{\pi^d 5^{n+1} d^{n/2}}
\ k^{-n/2}.
$$
\end{theorem}
This result holds for any set of vector generators.
On the other hand, for certain sets of {\em badly approximable}
generators (to be defined later), we can establish the following upper
bound.
\begin{theorem}
\label{upperbound}
Let $A_{n\times d}$ be a badly approximable matrix, with
rows $\vec\alpha_1,...,\vec\alpha_n$, and 
approximation constant $C_A$.  
If $Q$ is the generating measure of the 
random walk on the $d$-torus generated by the $\vec\alpha_i$,
then the $k$-th step
probability distribution $Q^{*k}$ satisfies:
$$
D(Q^{*k}) \le
        \left(\frac{3}{2}\right)^d 20
        \left(\frac{n}{C_A\sqrt{2}}\right)^{n/d} k^{-n/2d}.
$$
\end{theorem}
We note that the case $d=1$ corresponds to a random walk on the
circle, which has been studied for a single generator~\cite{su-tams}
and for several generators~\cite{hensley-su}.

\section{Lower Bound}
The following notation will be used throughout this paper:
\begin{description}
\item[$\|x\|$] the Euclidean ($L^2$) norm of a vector $x$
\item[$\|x\|_\infty$] the supremum norm of a vector $x$
\item[$\{x\}$] the Euclidean ($L^2$) distance 
from $x$ to the nearest integral point
\item[$\{x\}_\infty$] the supremum distance from $x$ to the nearest
  integral point
\end{description}

To establish a lower bound for the discrepancy,
we use a lemma due to Dirichlet:
\begin{lemma}[Dirichlet 1842]
\label{dirichlet-lemma}
Given any real $n \times d$ matrix $A$ and $q\ge 1$, there
is some ${\bf h} \in \Z^d$ such that $0 < \|{\bf h}\|_\infty \le q^{n/d}$
and $\{A{\bf h}\}_\infty < 1/q$.
\end{lemma}

A simple proof using a pigeonhole argument may be found in~\cite{schmidt}.
We now prove Theorem \ref{lowerbound}.

\begin{proof}
Su \cite{su-leveque} 
has shown that for any probability distribution $P$ on $\T^d$:
\begin{equation}
\label{su-torus-lowerbound}
D(P) \ge \sup_{{\bf r} \in (0,.5]^d}
\left[
\sum_{{\bf 0} \neq {\bf h} \in \Z^d}
|\hat{P}({\bf h})|^2 \prod_{i=1}^d
\left\{
\begin{array}{ll}
        \frac{\sin^2(2\pi h_i r_i)}{\pi^2h_i^2} & \textrm{if\ } h_i \neq 0 \\
        4r_i^2                                  & \textrm{if\ } h_i = 0
\end{array}
\right\}
\right]^{1/2}
\end{equation}
where $\hat{P}({\bf h})$ is the {\em Fourier transform} of 
$P$, i.e., $\hat{P}({\bf h})
=\int_{\T^d} e^{2\pi i {\bf h}\cdot{\bf x}} Q(d{\bf x})$.
We will use
this formula to bound $D(Q^{*k})$ where $Q$ is the generating measure
of our random walk.  Note that:
\begin{eqnarray*}
\hat{Q}({\bf h}) & = & \sum_{j=1}^n \frac{1}{2n}
                        (e^{2\pi i {\bf h} \cdot \vec\alpha_j} +
                        e^{-2\pi i {\bf h} \cdot \vec\alpha_j})         \\
                & = & \frac{1}{n} \sum_{j=1}^n
                        \cos(2\pi {\bf h} \cdot \vec\alpha_j).
\end{eqnarray*}

Since $\cos(2\pi x) = \cos(2\pi \{x\}) \ge 1-2\pi^2\{x\}^2$, we have
\begin{eqnarray*}
\hat{Q}({\bf h}) & \ge & \frac{1}{n} \sum_{j=1}^n
                        1-2\pi^2 \{\vec\alpha_j \cdot {\bf h}\}^2       \\
                & \ge & 1-\frac{2\pi^2}{n}\sum_{j=1}^n
                                \{\vec\alpha_j \cdot {\bf h}\}^2        \\
                & \ge & 1-\frac{2\pi^2}{n} \{A {\bf h}\} ^2
\end{eqnarray*}
where $A\equiv(\vec\alpha_1\vec\alpha_2\cdots\vec\alpha_n)^\top$
is the $n \times d$ matrix whose rows are the $\alpha_j$'s. Also, noting that
$\widehat{Q^{*k}} ({\bf h}) = \hat{Q}^k ({\bf h})$ and that
$(1-x)^k \ge 1-kx$ for $k \ge 1$ and $x \le 1$, we have that
$$
\widehat{Q^{*k}} ({\bf h}) \ge 1-\frac{2\pi^2k}{n}\{A{\bf h}\}^2
$$
as long as $2\pi^2k\{A{\bf h}\}^2/n < 1$.
This is ensured by setting $Z_1 = 2\pi^2/25 < 1$ and letting $q =
(2\pi^2kd/Z_1)^{1/2}$. Then Lemma~\ref{dirichlet-lemma}
implies that there exists ${\bf h} \in \Z^d$ such that
$0 < \|{\bf h}\|_\infty \le q^{n/d}$ and $\{A{\bf h}\}_\infty < 1/q$.
This yields $2\pi^2k\{A{\bf h}\}^2/n \le 2\pi^2kd\{A{\bf h}\}^2_\infty
< 2\pi^2kd/q^2 = Z_1 < 1$, as desired.
(Note that $|\hat{Q}({\bf h})|^k \ge 1-Z_1$.) By evaluating inequality 
(\ref{su-torus-lowerbound}) at this
{\bf h} we find
$$
D(Q^{*k}) \ge \sup_{{\bf r} \in (0,.5]^d} \left[
(1-Z_1)\prod_{i=1}^d \left\{
\begin{array}{ll}
        \frac{\sin(2\pi h_ir_i)}{\pi h_i} & \textrm{if\ } h_i \neq 0 \\
        2r_i                              & \textrm{if\ } h_i = 0
\end{array}
\right\} \right]
$$

Then, if we let
$r_i = 1/4h_i < 1/2$ if $h_i \neq 0$ and $r_i = 1/2\pi$ if $h_i = 0$
and define $R({\bf h})=\prod_{i=1}^n \max\{1,|h_i|\}$ to relate the
size of ${\bf h}$,
we find that
\begin{eqnarray*}
D(Q^{*k})
  & \ge & (1-Z_1) \prod_{i=1}^d \left\{
          \begin{array}{ll}
                \frac{1}{\pi h_i} & \textrm{if\ } h_i \neq 0 \\
                \frac{1}{\pi}     & \textrm{if\ } h_i = 0
          \end{array}
          \right\}                                                      \\
  & \ge & \frac{1-Z_1}{\pi^d R({\bf h})}                                \\
  & \ge & \frac{1-Z_1}{\pi^d \|{\bf h}\|_\infty^d}                      \\
  & \ge & \frac{1-Z_1}{\pi^d q^n}                                       \\
  & \ge & \frac{1-Z_1}{\pi^d (2\pi^2kd/Z_1)^{n/2}}                      \\
  & \ge & \frac{1-Z_1}{\pi^d 5^n d^{n/2}} \ k^{-n/2}       \\
  & \ge & \frac{1}{\pi^d 5^{n+1} d^{n/2}} \ k^{-n/2}.
\end{eqnarray*}
\end{proof}

\section{Upper Bound}

We now seek an upper bound on the discrepancy of the random walk when
our generators arise as rows of a {\em badly approximable matrix}.

\begin{defn}
\label{bad-approx-defn}
We say an $n \times d$ matrix $A$ is {\em badly approximable} if
there exists a constant $C_A$ such that
$\{A{\bf h}\}_\infty > C_A/\|{\bf h}\|_\infty^{d/n}$ for all non-zero
${\bf h} \in \Z^d$.  We call $C_A$ the {\em approximation constant} of $A$.
\end{defn}

Note that Lemma~\ref{dirichlet-lemma} implies that for any matrix $A$ (not
just badly approximable ones),
$\{A{\bf h}\}_\infty < 1/\|{\bf h}\|_\infty^{d/n}$ for infinitely many 
${\bf h} \in \Z^d$.  
Thus we say $A$ is badly approximable if the reverse inequality
holds (up to a constant $C_A$) for all ${\bf h} \in \Z^d$.
This definition closely follows Schmidt \cite{schmidt69}, who
defines {\em badly approximable linear forms}; this 
corresponds to our definition by noting $A{\bf h}$ is a linear form in
the variables $h_i$.

As a subset of $\R^{nd}$, the set of badly approximable matrices has
Lebesgue measure zero~\cite{khintchine} although their Hausdorff
dimension is $nd$ and there are uncountably many of them~\cite{schmidt69}.

We now prove Theorem \ref{upperbound}.

\begin{proof}
It is known~\cite{drmota-tichy} from Erd\H{o}s, Tur\`an, and Koksma that
for all positive integers $M$,
\begin{eqnarray}\label{ineq}
D(Q^{*k}) \le
\left( \frac{3}{2} \right)^d
\left( \frac{2}{M+1} +
\sum_{{\bf h}\in\Z^d \atop 0 < \|{\bf h}\|_\infty \le M} 
        \frac{|\hat{Q}^k({\bf h})|}{R({\bf h})}
\right).
\end{eqnarray}

Since $|\cos(2\pi x)| \le 1-4\{2x\}^2$ for all $x \in \R$, it follows that
\begin{eqnarray*}
        |\hat{Q}({\bf h})|
                & = & \frac{1}{n} \sum_{j=1}^n | \cos(2\pi 
                                {\bf h} \cdot \vec\alpha_j) |   \\
                & \le & \frac{1}{n} \sum_{j=1}^n 1 - 4 \{ 2
                                {\bf h} \cdot \vec\alpha_j \}^2 \\
                & \le & 1-\frac{4}{n}\sum_{j=1}^n \{ 2
                                \vec\alpha_j \cdot {\bf h} \}^2 \\
                & \le & 1-\frac{4}{n} \{ 2 A {\bf h} \}^2       \\
                & \le & \exp\left(-\frac{4}{n} \{ 2 A {\bf h} \}^2\right).
\end{eqnarray*}

In light of
inequality~(\ref{ineq}),
we need to estimate a sum of the form
$$
\sum_{0 < \|{\bf h}\|_\infty \le M}
        \frac{|\hat{Q}^k({\bf h})|}{R({\bf h})}
\le \sum_{{\bf h}\in\Z^d \atop 0 < \|{\bf h}\|_\infty \le M}
        \frac{\exp\left(-\frac{4k}{n} \{2A{\bf h}\}^2\right)}
        {R({\bf h})}=:S.
$$

Since $M$ may be chosen freely, choose an integer $M$ such that
\begin{equation}
\label{M-defn}
M \le \frac{1}{8} \left(\frac{2k(C_A)^2}{n^2}\right)^{n/2d} < M+1.
\end{equation}
Here $C_A$ is an approximation constant for the badly approximable $A$
and $k$ is the number of steps in the walk.  We can show:
\begin{lemma}
\label{S-bound}
With $S$ and $M$ defined as above,
$S \leq \frac{0.5}{M+1}$.
\end{lemma}

Before proving this lemma, we show how the theorem follows.
From 
inequality~(\ref{ineq}),
we find
$$
D(Q^{*k}) \le \left(\frac{3}{2}\right)^d\left(\frac{2}{M+1}+ S\right) \le
\left(\frac{3}{2}\right)^d\left(\frac{2}{M+1}+\frac{0.5}{M+1}\right) =
\left(\frac{3}{2}\right)^d\frac{2.5}{M+1}.
$$
From the inequality for $M+1$ in (\ref{M-defn}), we have 
$$
D(Q^{*k}) \le \left(\frac{3}{2}\right)^d 20
        \left(\frac{n^2}{2k(C_A)^2}\right)^{n/2d}
 = \left(\frac{3}{2}\right)^d 20
        \left(\frac{n}{C_A\sqrt{2}}\right)^{n/d} k^{-n/2d},
$$
which concludes the proof of Theorem \ref{upperbound}.
\end{proof}

All that remains is to prove Lemma \ref{S-bound}.

\begin{proof}[Proof of Lemma \ref{S-bound}]
We shall bound $S$ in three stages: (1) first, we group the terms of $S$
into ``cohorts'' based on the size of $\|{\bf h}\|_\infty$,
(2) we note that the points $2A{\bf h}$ are bounded away from each
other in $\T^d$ and therefore can 
bound the terms within each cohort based on the size of $\{2A{\bf h}\}_\infty$,
and (3) estimating the resulting expression.

\subsubsection*{(1) 
Grouping the terms of $S$ by the size of $\|{\bf h}\|_\infty$}

Choose an integer $J$ such that
$$
2^{J-1} \le M \le 2^J-1.
$$
The sum in $S$ may be grouped into $J$ cohorts of integers
$\|{\bf h}\|_\infty \in H_j := [2^{j-1},2^j-1]$
for $j=\{1,\ldots,J\}$. Therefore,
\begin{eqnarray*}
S & \le & \sum_{j=1}^J
        \sum_{\|{\bf h}\|_\infty\in H_j} 
        \frac{\exp(-\frac{4k}{n}\{2A{\bf h}\}^2)}{R({\bf h})} \\
  & \le & \sum_{j=1}^J
        \sum_{\|{\bf h}\|_\infty\in H_j} 
        \frac{\exp(-\frac{4k}{n}\{2A{\bf h}\}^2)}{2^{j-1}}    \\
  & \le & \sum_{j=1}^J \frac{1}{2^{j-1}}
        \sum_{\|{\bf h}\|_\infty\in H_j} 
        \exp\left(-\frac{4k}{n}\{2A{\bf h}\}^2\right)         \\
  & \le & \sum_{j=1}^J \frac{1}{2^{j-1}}
        \sum_{\|{\bf h}\|_\infty\in H_j} 
        \exp\left(-\frac{4k}{n}\{2A{\bf h}\}^2_\infty\right).
\end{eqnarray*}

\subsubsection*{(2) Bounding the terms within each cohort}

Within each cohort $[2^{j-1},2^j-1]$, 
since ${\bf h}$ is a non-zero integral vector, 
the use of Definition
\ref{bad-approx-defn} yields
$\{2A{\bf h}\}_\infty > C_A/\|2{\bf h}\|_\infty^{d/n}$ 
where $C_A$ is the approximation constant of the matrix $A$. 
Therefore,
each $2A{\bf h}$ is bounded away from any integral
point by $C_A/\|2 {\bf h}\|_\infty^{d/n}$. In fact, they are also bounded away
from each other, since if
$\|{\bf h}_1\|_\infty$, $\|{\bf h}_2\|_\infty \in [2^{j-1},2^j-1]$ 
and ${\bf h}_1\neq{\bf h}_2$, 
then
$\|{\bf h}_1 - {\bf h}_2\|_\infty \le 2^{j+1}$ and
$$
\{ 2A({\bf h}_1 - {\bf h}_2) \}_\infty >
\frac{C_A}{2^{d/n} \|{\bf h}_1 - {\bf h}_2\|_\infty^{d/n}} \ge
\frac{C_A}{2^{d/n} (2^{j+1})^{d/n}} = \frac{C_A}{2^{(j+2)d/n}}.
$$
Therefore, we divide the unit cube $[0,1]^n$ into subcubes of side-length
$C_A/(2^{d(j+2)/n})$ and distribute the points $\{2A{\bf h}\}$
throughout them. In the worst case, all of the points are
distributed near the corners of the cube and occupy adjacent subcubes.
Therefore,
\begin{eqnarray*}
S & \le & \sum_{j=1}^J \frac{2^n}{2^{j-1}} \sum_{i=1}^{|H_j|^{1/n}}
        (i+1)^{n-1} \exp \left(-\frac{4k}{n}
        \left(\frac{i C_A}{2^{(j+2)d/n}} \right)^2 \right) \\
  & \le & \sum_{j=1}^J 2^{n+1-j} \sum_{i=1}^\infty
        (i+1)^{n-1} \exp\left(-\frac{4 k i^2 (C_A)^2}{n 2^{2(j+2)d/n}}\right).
\end{eqnarray*}

\subsubsection*{(3) Estimating the resulting expression.}

Since $M \ge 2^{J-1}$ and
$$
k \ge \frac{n^2 2^{6d/n}M^{2d/n}}{2 (C_A)^2}
  \ge \frac{n^2 2^{6d/n}(2^{J-1})^{2d/n}}{2 (C_A)^2}
  = \frac{n^2 2^{2(J+2)d/n}}{2 (C_A)^2},
$$
we can say that:
$$
S \le \sum_{j=1}^J 2^{n+1-j} \sum_{i=1}^\infty
        (i+1)^{n-1} \exp\left(-2 i^2 n (2^{J-j})^{2d/n}\right).
$$

Since $j \le J$, $i \ge 1$ and $n \ge 1$, the log derivative with
respect to $i$ of the inner sum can be bounded:
$$
\frac{n-1}{i+1} - 4i n (2^{J-j})^{2d/n} \le \frac{n-1}{i+1} - 4 i n \le
-4 i \le -4.
$$

Therefore, the expression in the inner sum decreases geometrically by
at least the ratio $e^{-4}$, and the inner sum can be bounded by the first
term (at $i = 1$) times the constant $1/(1-e^{-4})$. 
Therefore,
\begin{eqnarray*}
S & \le & \sum_{j=1}^J \frac{2^{n+1-j}}{1-e^{-4}}
                2^{n-1}\exp\left(-2n2^{2(J-j)d/n}\right) \\
  & \le & \frac{2^{2n}}{1-e^{-4}}\sum_{j=1}^J 2^{-j}
                \exp\left(-2n 2^{2(J-j)d/n}\right).
\end{eqnarray*}

This sum may be bounded by noting that the largest term occurs when $j=J$.
For $j \le J$, the log derivative of the terms with respect to $j$ is
$\ln 2(-1+d2^{2+2(J-j)d/n}) \ge \ln 2 (-1+ 2^2)
= \ln 8$. Therefore, the sum decreases geometrically
with ratio at least $1/8$ as $j \le J$ decreases, so the sum is bounded
by seven-eighths the final term at $j=J$. Also, recalling that $M \le 2^J-1$,
\begin{eqnarray*}
S & \le & \frac{2^{2n}}{1-e^{-4}} \frac{7}{8} 2^{-J}
        \exp\left(-2n2^{2(J-J)d/n}\right) =
\frac{7}{8(1-e^{-4})}\left(\frac{2}{e}\right)^{2n}\frac{1}{M+1} \\
  & \le & \frac{28}{8(e^2-e^{-2})}\frac{1}{M+1} \le
\frac{0.5}{M+1},
\end{eqnarray*}
as was to be shown.
\end{proof}

So, for badly approximable matrices $A$, we have the following
discrepancy bounds on the associated random walk:
$$
\frac{C(n,d)}{k^{n/2}} \le D(Q^{*k}) \le \frac{C(n,d,A)}{k^{n/2d}}.
$$
In general, these bounds do not match unless $d=1$.  In that case, we
recover the same order of convergence as in \cite{hensley-su}.  

We conjecture that 
in the lower bound of Theorem
\ref{lowerbound} (which applies to any random walk on the torus generated
by a finite set of vectors) the $k^{n/2}$ 
may be improved to $k^{n/2d}$ (for all matrices, not just badly
approximable ones).
This would yield matching upper and
lower bounds for random walks on the $d$-torus generated by the rows of a 
badly approximable matrix.  This would confirm that such walks converge
the fastest among all finitely-generated random walks on the $d$-torus,
a fact that has already been shown in \cite{hensley-su} for dimension $d=1$.
Developing an approximation lemma similar to Dirichlet's Lemma 
(Lemma \ref{dirichlet-lemma}) 
that bounds $R({\bf h})$ instead of $\|{\bf h}\|_\infty$ may help in
this regard.

\end{document}